\documentclass[a4paper,10pt]{article}
\usepackage[latin1]{inputenc}
\usepackage[T1]{fontenc}
\usepackage{amsmath}
\usepackage{amsfonts}
\usepackage{amstext}
\usepackage{amsthm}
\usepackage[all]{xy}
\usepackage{geometry}
\usepackage{srcltx}
\usepackage{graphics}
\usepackage{epsfig}
\usepackage{subfigure}
\usepackage{slashed}
\usepackage{color}
\usepackage{amssymb}
\usepackage{srcltx}
\newtheorem{theorem}{Theorem}[section]

\newtheorem{cor}[theorem]{Corollary}

\DeclareMathOperator{\grad}{grad}

\DeclareMathOperator{\sgn}{sgn}

\DeclareMathOperator{\diverg}{div}
\DeclareMathOperator{\dvol}{dvol}

\title{On bifurcation for semilinear elliptic Dirichlet problems on geodesic balls}
\author{Alessandro Portaluri and Nils Waterstraat}

\begin{document}
\date{}
\maketitle

\footnotetext[2]{A. Portaluri was supported by the grant PRIN2009 ``Critical Point Theory and Perturbative Methods
for Nonlinear Differential Equations''.}
\footnotetext[3]{N. Waterstraat was supported by a postdoctoral fellowship of
the German Academic Exchange Service (DAAD).}

\begin{abstract}
We study bifurcation from a branch of trivial solutions of semilinear elliptic Dirichlet boundary value problems on a geodesic ball, whose radius is used as the bifurcation parameter. In the proof of our main theorem we obtain in addition a special case of an index theorem due to S. Smale.
\end{abstract}

\section{Introduction}
Let $(M,g)$ be an oriented Riemannian manifold of dimension $n$ and let $\Delta=\diverg\grad:C^\infty(M)\rightarrow C^\infty(M)$ denote the associated Laplace-Beltrami operator. Let $V:M\times\mathbb{R}\rightarrow\mathbb{R}$ be a smooth function such that $V(p,0)=0$ for all $p\in M$ and 

\begin{align}\label{growth}
|V(p,\xi)|\leq C(1+|\xi|^\alpha),\quad \left|\frac{\partial V}{\partial\xi}(p,\xi)\right|\leq C(1+|\xi|^\beta),\quad (p,\xi)\in M\times\mathbb{R}, 
\end{align}
for some $C>0$ and constants $\alpha,\beta\geq 0$ depending on the dimension $n$ of $M$ (cf. \cite[\S 1.2]{Ambrosetti}). In this paper we deal with local solutions of the semilinear equation

\begin{align}\label{equation}
-\Delta u(p)+V(p,u(p))=0,\quad p\in M,
\end{align}
under Dirichlet boundary conditions. Note that many equations from geometric analysis are of the type \eqref{equation}. Let us refer to \cite{Aubin}, \cite{Besse} and just mention as an example on compact manifolds of dimension $n\geq 3$ the equation

\begin{align}\label{equationcurvature}
4\,\frac{n-1}{n-2}\,\Delta u(p)+s(p)u(p)=\mu\, u(p)^\frac{n+2}{n-2},\quad p\in M,
\end{align}
where $s:M\rightarrow\mathbb{R}$ denotes the scalar curvature function and $\mu$ the Yamabe invariant of the metric $g$ on $M$. Positive solutions $u\in C^\infty(M)$ of \eqref{equationcurvature} give rise to metrics $\tilde{g}$ of constant scalar curvature on $M$ by $\tilde{g}=u^\frac{4}{n-2}g$.\\
We now fix a point $p_0\in M$ and assume that the unit ball $B(0,1)\subset T_{p_0}M$ is contained in the maximal domain on which the exponential map $\exp_{p_0}$ at $p_0$ is an embedding. Let us denote by $B(p_0,r)=\exp_{p_0}(B(0,r))$ the geodesic ball of radius $0<r\leq 1$ around $p_0$ and consider the Dirichlet boundary value problems

\begin{equation}\label{bvp}
\left\{
\begin{aligned}
-\Delta u(p)+V(p,u(p))&=0,\quad p\in B(p_0,r)\\
u(p)&=0,\quad p\in\partial B(p_0,r).
\end{aligned}
\right.
\end{equation}
We call $r^\ast\in(0,1]$ a \textit{bifurcation radius} for the boundary value problems \eqref{bvp} if 
there exists a sequence of radii $r_n\rightarrow r^\ast$ and functions $u_n\in H^1_0(B(p_0,r_n))$ such that $u_n$ is a non-trivial weak solution of \eqref{bvp} on $B(p_0,r_n)$ and $\|u_n\|_{H^1_0(B(p_0,r_n))}\rightarrow 0$. Note that we exclude from the definition the limiting case $r^\ast=0$ in which the domain degenerates to a point. The reason is that $\|u_n\|_{H^1_0(B(p_0,r_n))}\rightarrow 0$ for $r_n\rightarrow 0$ holds, for example, for any sequence of functions $u_n\in C^1(B(p_0,r_n))$, $n\in\mathbb{N}$, such that all $u_n$ and their weak derivatives are bounded uniformly. Consequently, a bifurcation radius $r^\ast=0$ would not imply the existence of non-trivial solutions of \eqref{bvp} for small $r>0$ which are arbitrarily close to the trivial solution $u\equiv 0$ in a suitable sense.\\
Let us now consider the linearised boundary value problems

\begin{equation}\label{bvpII}
\left\{
\begin{aligned}
-\Delta u(p)+f(p)u(p)&=0,\quad p\in B(p_0,r)\\
u(p)&=0,\quad p\in\partial B(p_0,r),
\end{aligned}
\right.
\end{equation}
where $f(p)=\frac{\partial V}{\partial \xi}(p,0)$, $p\in M$. We call $r^\ast\in(0,1]$ a \textit{conjugate radius} for \eqref{bvpII} if

\begin{align*}
m(r^\ast):=\dim\{u\in C^2(B(p_0,r^\ast)):\, u\,\,\,\text{solves}\,\,\,\eqref{bvpII}\}>0,
\end{align*}
and from now on we assume that $m(1)=0$. Our main result reads as follows:

\begin{theorem}\label{theorem}
The bifurcation radii of \eqref{bvp} are precisely the conjugate radii of \eqref{bvpII}. 
\end{theorem}

We explain below in the proof of Theorem \ref{theorem} that we obtain from our methods a new proof of the Morse-Smale index theorem \cite{Smale} (cf. also \cite{SmaleCorr}) for the linearised equations \eqref{bvpII}. As a consequence, we conclude that $m(r)=0$ for almost all radii $r\in(0,1)$, and moreover, we derive from Theorem \ref{theorem} the following corollary:

\begin{cor}\label{cor}
Let $\mu$ denote the Morse index of \eqref{bvpII} on $B(p_0,1)$, i.e. the number of negative eigenvalues counted according to their multiplicities. If $\mu\neq 0$, then there exist at least
	
	\begin{align*}
	\left\lfloor\frac{\mu}{\max_{0<r<1}m(r)}\right\rfloor
	\end{align*}
distinct bifurcation radii in $(0,1)$, where $\lfloor\cdot\rfloor$ denotes the integral part of a real number.
\end{cor}
Let us point out that a proof of Theorem \ref{theorem} and Corollary \ref{cor} for the special case that $M$ is a star-shaped domain in $\mathbb{R}^n$ can be found in \cite{AleIchDomain}. The following section is devoted to the more general setting which we consider here.


\section{The proof}
Our main reference for the Laplace-Beltrami operator on manifolds with boundary is \cite[\S 2.4]{Taylor}. Let us recall at first that in local coordinates

\[\Delta u=\sum^n_{j,k=1}{|g|^{-\frac{1}{2}}\frac{\partial}{\partial x^j}\left(g^{jk}|g|^\frac{1}{2}\frac{\partial u}{\partial x^k}\right)},\]
where $g^{jk}$, $1\leq j,k\leq n$, are the components of the inverse of the metric tensor $g=\{g_{jk}\}$ and $|g|:=|\det\{g_{jk}\}|$ is the absolute value of its determinant. Denoting by $\dvol_g$ the volume form of $g$, we find for $v\in H^1_0(B(p_0,r))$, $0<r\leq 1$,

\begin{align*}
&-\int_{B(p_0,r)}{(\Delta u)(p)v(p)\,\dvol_g}+\int_{B(p_0,r)}{V(p,u(p))\,v(p)\,\dvol_g}\\
&=-\int_{B(0,r)}{v(x)\sum^n_{j,k=1}{\frac{\partial}{\partial x^j}\left(g^{jk}(x)|g(x)|^\frac{1}{2}\frac{\partial u}{\partial x^k}(x)\right)}dx}+\int_{B(0,r)}{|g(x)|^\frac{1}{2}V(x,u(x))v(x)\, dx}\\
&=\int_{B(0,r)}{\sum^n_{j,k=1}{g^{jk}(x)|g(x)|^\frac{1}{2}\frac{\partial u}{\partial x^k}(x)\frac{\partial v}{\partial x^j}(x)}\,dx}+\int_{B(0,r)}{|g(x)|^\frac{1}{2}V(x,u(x))v(x)\, dx}\\
&=r\int_{B(0,1)}{\sum^n_{j,k=1}{g^{jk}(r\cdot x)|g(r\cdot x)|^\frac{1}{2}\frac{\partial u}{\partial x^k}(r\cdot x)\frac{\partial v}{\partial x^j}(r\cdot x)}\,dx}\\
&+r\int_{B(0,1)}{|g(r\cdot x)|^\frac{1}{2}V(r\cdot x,u(r\cdot x))v(r\cdot x)\, dx},
\end{align*}
and analogously

\begin{align*}
&-\int_{B(p_0,r)}{(\Delta u)(p)v(p)\,\dvol_g}+\int_{B(p_0,r)}{f(p)u(p)v(p)\,\dvol_g}\\
&=r\int_{B(0,1)}{\sum^n_{j,k=1}{g^{jk}(r\cdot x)|g(r\cdot x)|^\frac{1}{2}\frac{\partial u}{\partial x^k}(r\cdot x)\frac{\partial v}{\partial x^j}(r\cdot x)}\,dx}\\
&+r\int_{B(0,1)}{|g(r\cdot x)|^\frac{1}{2}f(r\cdot x)u(r\cdot x)v(r\cdot x)\, dx}.
\end{align*}
We now set $B:=B(0,1)$ and define for $0\leq r\leq 1$ a functional $q_r:H^1_0(B)\times H^1_0(B)\rightarrow\mathbb{R}$ by

\[q_r(u,v)=\int_{B}{\sum^n_{j,k=1}{g^{jk}(r\cdot x)|g(r\cdot x)|^\frac{1}{2}\frac{\partial u}{\partial x^k}(x)\frac{\partial v}{\partial x^j}(x)}\,dx}+r^2\int_{B}{|g(r\cdot x)|^\frac{1}{2}V(r\cdot x,u(x))v(x)\, dx}\]
as well as a quadratic form $h_r:H^1_0(B)\rightarrow\mathbb{R}$ by

\[h_r(u)=\int_{B}{\sum^n_{j,k=1}{g^{jk}(r\cdot x)|g(r\cdot x)|^\frac{1}{2}\frac{\partial u}{\partial x^k}(x)\frac{\partial u}{\partial x^j}(x)}\,dx}+r^2\int_{B}{|g(r\cdot x)|^\frac{1}{2}f(r\cdot x)u(x)^2\, dx}.\]
From the computations above we conclude that:

\begin{enumerate}
	\item[i)] $r^\ast\in(0,1]$ is a bifurcation radius for \eqref{bvp}, if and only if there exist a sequence $\{r_n\}_{n\in\mathbb{N}}\subset(0,1]$, $r_n\rightarrow r^\ast$, and a sequence of non-trivial functions $\{u_n\}_{n\in\mathbb{N}}\subset H^1_0(B)$, $u_n\rightarrow 0$, such that $q_{r_n}(u_n,\cdot)=0\in (H^1_0(B))^\ast$ for all $n\in\mathbb{N}$.
	\item[ii)] $r^\ast\in(0,1]$ is a conjugate radius for \eqref{bvpII} if and only if $h_r$ is degenerate.
\end{enumerate} 
We now define a function $\psi:[0,1]\times H^1_0(B)\rightarrow\mathbb{R}$ by

\[\psi(r,u)=\int_{B}{\sum^n_{j,k=1}{g^{jk}(r\cdot x)|g(r\cdot x)|^\frac{1}{2}\frac{\partial u}{\partial x^k}(x)\frac{\partial u}{\partial x^j}(x)}\,dx}+r^2\int_{B}{|g(r\cdot x)|^\frac{1}{2}\,G(r\cdot x,u(x))\,dx},\] 
where

\[G(x,t)=\int^t_0{V(x,\xi)\,d\xi}.\]
It is a standard result that $\psi$ is $C^2$-smooth under the growth conditions \eqref{growth}, and $D_u\psi_r=q_r(u,\cdot)$, $u\in H^1_0(B)$. Moreover, $0\in H^1_0(B)$ is a critical point of all functionals $\psi_r$ and the corresponding Hessians are given by $D^2_0\psi_r=h_r$. From the compactness of the inclusion $H^1_0(B)\hookrightarrow L^2(B)$, we see at once that the Riesz representation of the quadratic form $h_r$ is a selfadjoint Fredholm operator. In particular, it is invertible if $h_r$ is non-degenerate.\\
Let us now assume at first that $r^\ast\in(0,1]$ is not a conjugate radius. Then $h_{r^\ast}$ is non-degenerate and we conclude by the implicit function theorem \cite[\S 2.2]{Ambrosetti} that the equation $q_r(u,\cdot)=0$ has no other solutions than $(r,0)\in [0,1]\times H^1_0(B)$ in a neighbourhood of $(r^\ast,0)$. Consequently, $(r^\ast,0)$ is not a bifurcation radius, and we have shown that every bifurcation radius in $(0,1]$ is a conjugate radius.\\
In order to prove the remaining implication of Theorem \ref{theorem}, we make use of the bifurcation theory for critical points of smooth functionals developed in \cite{SFLPejsachowicz}. Accordingly, we consider for $r_0\in(0,1)$ the quadratic forms

\[\Gamma(h,r_0):\ker h_{r_0}\rightarrow\mathbb{R},\quad \Gamma(h,r_0)[u]=\left(\frac{d}{dr}\mid_{r=r_0}h_r\right)u.\]   By \cite[Thm. 1\& Thm. 4.1]{SFLPejsachowicz}, $r_0$ is a bifurcation radius if $\Gamma(h,r_0)$ is non-degenerate and has a non-vanishing signature (cf. also Section 2.1 in \cite{AleIchDomain}). Consequently, we now assume that $r_0\in(0,1)$ is a conjugate radius and our aim is to compute $\Gamma(h,r_0)$. Let us write for simplicity of notation

\begin{align*}
a^{jk}(x)&=g^{jk}(x)|g(x)|^\frac{1}{2},\quad x\in B,\,\,1\leq j,k\leq n,\\
\tilde{f}(x)&=|g(x)|^\frac{1}{2}f(x),\quad x\in B.
\end{align*}
For $u\in\ker h_{r_0}$ we have by definition

\begin{align}\label{crossingform}
\Gamma(h,r_0)[u]=\int_B{\sum^n_{j,k=1}{\langle\nabla a^{jk}(r_0\cdot x),x\rangle\frac{\partial u}{\partial x^k}\frac{\partial u}{\partial x^j}\,dx}}+\int_B{\frac{d}{d r}\mid_{r=r_0}(r^2\tilde{f}(r\cdot x))u(x)^2\,dx}.
\end{align} 
Let us now introduce a new function by $v_r(x):=u(\frac{r}{r_0}\cdot x)$, $r\in(0,r_0]$, $x\in B$, and denote

\begin{align}\label{udot}
\dot{u}(x):=\frac{d}{dr}\mid_{r=r_0}v_r(x)=\frac{1}{r_0}\langle\nabla u(x),x\rangle.
\end{align}
It is readily seen that $v_r$ satisfies

\begin{align*}
-\sum^n_{j,k=1}{\frac{\partial}{\partial x^j}\left(a^{jk}(r\cdot x)\frac{\partial v_r}{\partial x^k}\right)}+r^2\tilde{f}(r\cdot x)v_r(x)=0,
\end{align*} 
and by differentiating with respect to $r$ at $r=r_0$ we have

\begin{align}\label{equ1}
\begin{split}
0&=-\sum^n_{j,k=1}{\frac{\partial}{\partial x^j}\left(\langle\nabla a^{jk}(r_0\cdot x),x\rangle\frac{\partial u}{\partial x^k}\right)}-\sum^n_{j,k=1}{\frac{\partial}{\partial x^j}\left(a^{jk}(r_0\cdot x)\frac{\partial\dot{u}}{\partial x^k}\right)}\\
&+\frac{d}{d r}\mid_{r=r_0}(r^2\tilde{f}(r\cdot x))u(x)+r^2_0\tilde{f}(r_0\cdot x)\dot{u}(x).
\end{split}
\end{align} 
We multiply \eqref{equ1} by $u$ and integrate over $B$:

\begin{align*}
0&=-\int_B{\sum^n_{j,k=1}{\frac{\partial}{\partial x^j}\left(\langle\nabla a^{jk}(r_0\cdot x),x\rangle\frac{\partial u}{\partial x^k}\right)}u(x)\,dx}-\int_B{\sum^n_{j,k=1}{\frac{\partial}{\partial x^j}\left(a^{jk}(r_0\cdot x)\frac{\partial\dot{u}}{\partial x^k}\right)}u(x)\,dx}\\
&+\int_B{\frac{d}{dr}\mid_{r=r_0}(r^2\tilde{f}(r\cdot x))u(x)^2\,dx}+\int_B{r^2_0\tilde{f}(r_0\cdot x)\dot{u}(x)u(x)\,dx}.
\end{align*} 
Let $\nu(x)=(\nu_1(x),\ldots,\nu_n(x))$, $x\in\partial B$, denote the outward pointing unit normal to the boundary of $B$. Using $u\mid_{\partial B}=0$, we obtain from integration by parts

\begin{align*}
0&=\int_B{\sum^n_{j,k=1}{\langle\nabla a^{jk}(r_0\cdot x),x\rangle\frac{\partial u}{\partial x^k}}\frac{\partial u}{\partial x^j}\,dx}-\int_{\partial B}{\left(\sum^n_{j,k=1}{\langle\nabla a^{jk}(r_0\dot x),x\rangle \nu_j(x)\frac{\partial u}{\partial x^k}}\right)u(x)\,dS}\\
&+\int_B{\sum^n_{j,k=1}{a^{jk}(r_0\cdot x)\frac{\partial\dot{u}}{\partial x^k}}\frac{\partial u}{\partial x^j}\,dx}-\int_{\partial B}{\left(\sum^n_{j,k=1}{a^{jk}(r_0\cdot x)\nu_j(x)\frac{\partial\dot{u}}{\partial x^k}}\right)u(x)\,dS}\\
&+\int_B{\frac{d}{dr}\mid_{r=r_0}(r^2\tilde{f}(r\cdot x))u(x)^2\,dx}+\int_B{r^2_0\tilde{f}(r_0\cdot x)\dot{u}(x)u(x)\,dx}\\
&=\int_B{\sum^n_{j,k=1}{\langle\nabla a^{jk}(r_0\cdot x),x\rangle\frac{\partial u}{\partial x^k}}\frac{\partial u}{\partial x^j}\,dx}-\int_B{\sum^n_{j,k=1}{\frac{\partial}{\partial x^j}\left(a^{jk}(r_0\cdot x)\frac{\partial u}{\partial x^k}\right)}\dot{u}(x)\,dx}\\
&+\int_{\partial B}{\left(\sum^n_{j,k=1}{a^{jk}(r_0\cdot x)\nu_j(x)\frac{\partial u}{\partial x^k}}\right)\dot{u}(x)\,dS}\\
&+\int_B{\frac{d}{dr}\mid_{r=r_0}(r^2\tilde{f}(r\cdot x))u(x)^2\,dx}+\int_B{r^2_0\tilde{f}(r_0\cdot x)\dot{u}(x)u(x)\,dx}.
\end{align*}
Since $u\in\ker h_{r_0}$,

\begin{align}\label{PDE}
-\sum^n_{j,k=1}{\frac{\partial}{\partial x^j}\left(a^{jk}(r_0\cdot x)\frac{\partial u}{\partial x^k}\right)}
+r^2_0\tilde{f}(r_0\cdot x)u(x)=0,\quad x\in B,
\end{align}
and it follows from \eqref{crossingform} and \eqref{udot} that

\begin{align}\label{crossingformII}
\Gamma(h,r_0)[u]&=-\frac{1}{r_0}\int_{\partial B}{\left(\sum^n_{j,k=1}{a^{jk}(r_0\cdot x)\nu_j(x)\frac{\partial u}{\partial x^k}}\right)\langle \nabla u(x),x\rangle\,dS}.
\end{align}
If we set $A(x):=\{a_{jk}(x)\}$, $x\in B$, and use that $\nu(x)=x$ for all $x\in\partial B$, we can rewrite \eqref{crossingformII} as

\begin{align*}
\Gamma(h,r_0)[u]=-\frac{1}{r_0}\int_{\partial B}{\langle A(r_0\cdot x)x,\nabla u(x) \rangle\, \langle \nabla u(x),x\rangle\,dS}.
\end{align*}
Denoting by $(A(r_0\cdot x)x)^T$, $x\in\partial B$, the tangential component of the vector $A(r_0\cdot x)x$, we have

\[\langle A(r_0\cdot x)x,\nabla u(x)\rangle=\langle\nabla u(x),x\rangle\,\langle A(r_0\cdot x)x,x\rangle+\langle\nabla u(x),(A(r_0\cdot x)x)^T\rangle\] 
and hence

\begin{align*}
\Gamma(h,r_0)[u]&=-\frac{1}{r_0}\int_{\partial B}{\langle\nabla u(x),x\rangle^2\,\langle A(r_0\cdot x)x,x \rangle\,dS}\\
&-\frac{1}{r_0}\int_{\partial B}{\langle\nabla u(x),x\rangle\,\langle\nabla u(x), (A(r_0\cdot x)x)^T \rangle\,dS}.
\end{align*}
Since

\[\langle\nabla u(x),x\rangle\,\langle\nabla u(x),(A(r_0\cdot x)x)^T \rangle=\diverg(u(x)\langle x,\nabla u(x)\rangle(A(r_0\cdot x)x)^T),\quad x\in\partial B,\]
we finally get by Stokes' theorem

\begin{align}\label{crossfinal}
\Gamma(h,r_0)[u]=-\frac{1}{r_0}\int_{\partial B}{\langle\nabla u(x),x\rangle^2\,\langle A(r_0\cdot x)x,x \rangle\,dS}\leq 0,
\end{align}
where we use that $A(x)$ is positive definite for all $x\in B$.\\
Moreover, we obtain from \eqref{crossfinal} that if $\Gamma(h,r_0)[u]=0$ for some $u\in\ker h_{r_0}$, then 

\[\langle\nabla u(x),x\rangle=\langle\nabla u(x),\nu(x)\rangle=\frac{\partial u}{\partial\nu}(x)=0\]
for all $x\in\partial B$ which implies $u\equiv 0$ by the unique continuation property.\\
In summary, we have shown that $\Gamma(h,r_0)$ is negative definite, and so in particular non-degenerate with the non-vanishing signature

\begin{align}\label{crossingformIII}
\sgn\Gamma(h,r_0)=m(r_0).
\end{align}
Consequently, $r_0$ is a bifurcation radius and Theorem \ref{theorem} is proven.\\
Let us now prove Corollary \ref{cor}. We note at first that the Morse index $\mu$ of \eqref{bvpII} on the full domain $B(p_0,1)$ is given by the Morse index $\mu(h_1)$ of the quadratic form $h_1$. Moreover, since $h_0$ is positive, we see that $\mu(h_0)=0$. It is shown in \cite[Prop. 3.9\& Thm. 4.1]{SFLPejsachowicz} that if $\Gamma(h,r)$ is non-degenerate for all $r\in(0,1)$, then $\ker h_r=0$ for almost all $r\in(0,1)$ and

\[\mu(h_1)-\mu(h_0)=\sum_{0<r<1}{\sgn\Gamma(h,r)}.\]
Consequently, we conclude from \eqref{crossingformIII} that $m(r)=\dim\ker h_r=0$ for almost all $r\in(0,1)$ and

\begin{align}\label{Smale}
\mu=\sum_{0<r<1}{m(r)}.
\end{align}
Let us point out that \eqref{Smale} was obtained by Smale in \cite{Smale} by studying the monotonicity of eigenvalues under shrinking of domains. Hence we have obtained a new proof of Smale's theorem for the boundary value problem \eqref{bvpII}, and moreover, Corollary \ref{cor} is now an immediate consequence of \eqref{Smale} and Theorem \ref{theorem}.

\vspace{1cm}
Alessandro Portaluri\\
Department of Agriculture, Forest and Food Sciences\\
Universit\`a degli studi di Torino\\
Via Leonardo da Vinci, 44\\
10095 Grugliasco (TO)\\
Italy\\
E-mail: alessandro.portaluri@unito.it

\vspace{1cm}
Nils Waterstraat\\
Dipartimento di Scienze Matematiche\\
Politecnico di Torino\\
Corso Duca degli Abruzzi, 24\\
10129 Torino\\
Italy\\
E-mail: waterstraat@daad-alumni.de

\end{document}